\documentclass[11pt]{article}
\usepackage{url}
\usepackage{amssymb}
\usepackage{amsthm}
\usepackage{amsmath}

\newcommand{\dT}{d_{\mathrm{T}}}
\newcommand{\DD}{\mathfrak{D}}
\newcommand{\Fm}{\mathcal{F}_{(m)}}
\newcommand{\hook}{\kern 3pt \vrule height 0pt depth 0.4pt width 3pt
   \vrule height 5pt depth 0.4pt\kern 3pt}

\newcommand{\iT}{i_{\mathrm{T}}}

\newcommand{\Obar}{\overline{\Omega}}
\newcommand{\pd}[2]{\frac{\partial #1}{\partial #2}}
\newcommand{\R}{\mathbb{R}}
\newcommand{\St}{\widetilde{S}}
\newcommand{\vf}[1]{\frac{\partial}{\partial #1}}
\newcommand{\X}{\mathfrak{X}}

\newcommand{\ce}{\mathcal{E}}
\renewcommand{\O}{\Omega}

\newcommand{\art}[6]{#1: #2 {\it #3\/} {\bf #4} (#5) #6}

\newcommand{\inbook}[6]{#1: #2 {\it In:\ #3\/} (#4, #5) #6}

\begin{document}

\title{The fundamental form of a homogeneous Lagrangian
in two independent variables}
\author{D. J. Saunders\\
Department of Algebra and Geometry, 
Palacky Unversity \\
779 00 Olomouc, Czech Republic\footnote{Address for
correspondence: 30 Little Horwood Road, Great Horwood,
Milton Keynes, MK17 0QE, UK}\\
e-mail \url{david@symplectic.demon.co.uk}\\[2ex]
and\\[2ex]
M.\ Crampin\\
Department of Mathematical Physics and Astronomy\\
Ghent University, Krijgslaan 281, B--9000 Gent, Belgium\\
and\\
Department of Mathematics, King's College\\
Strand, London WC2R 2LS, UK}

\date{}

\maketitle

\begin{abstract}\noindent
We construct, for a homogeneous Lagrangian of arbitrary order in two independent
variables, a differential 2-form with the property that it is closed precisely
when the Lagrangian is null. This is similar to the property of the
`fundamental Lepage equivalent' associated with first-order Lagrangians
defined on jets of sections of a fibred manifold.
\\[1ex]
{\bf Keywords:} homogeneous Lagrangian, Lepage equivalent, Euler-Lagrange
form\\[1ex]
{\bf MSC2000 Classification:} 58E99, 49F99
\end{abstract}

\section{Introduction}

The `Lepage equivalents' of a Lagrangian are important tools for use when
studying variational problems on fibred manifolds:\ they are differential forms
having the same extremals as the Lagrangian form, with a further property
ensuring that their differentials give rise to the Euler-Lagrange form. The
classical example of a Lepage equivalent is the Cartan form in mechanics. In this
paper we study differential forms which play a r\^{o}le similar to that of Lepage
forms but in the context of homogeneous variational problems; there is a close
relationship, described below, between the two types of problem, and we believe
that studying the homogeneous context can throw some light on an existing unsolved problem
regarding Lepage equivalents of null Lagrangians.

In the context of a fibred manifold $\pi : E \to M$ with $\dim M = m$, a Lagrangian
is an $m$-form $\lambda \in \Omega^m J^k\pi$; any Lepage equivalent
$\theta$ of $\lambda$ will be defined on a jet manifold $J^l\pi$ (with, in general,
$l \geq k$) and will satisfy the conditions that $\theta - \pi_{l,k}^* \lambda$
should be contact, and that for any vector field $Z \in \X(J^l\pi)$ vertical
over $E$ the contraction $i_Z d\theta$ should also be contact. The
Euler-Lagrange form $\varepsilon$ is then the 1-contact part of $d\theta$. In mechanics,
for example, if we take coordinates $(t, q^a, \dot{q}^a)$ and a Lagrangian
$\lambda = L\, dt$ then
\[
\theta = L \, dt + \pd{L}{\dot{q}^a} (dq^a - \dot{q}^a dt) \, .
\]
See, for instance,~\cite{Gar,GS,Got,Krup2} for various approaches to the construction
of Lepage equivalents.

Global Lepage equivalents may always be found for a given Lagrangian, and if
$m = 1$ then they are unique. They are never unique when $m > 1$, because
adding an arbitrary non-zero 2-contact form to any Lepage equivalent will give
a different Lepage equivalent, although such a modification will not affect the
Euler-Lagrange form. Nevertheless, when the order $k$ of the Lagrangian is no
more than 2 then it is possible to make a canonical choice of Lepage
equivalent; this cannot, however, be done when $k \geq 3$ without the
specification of some additional structure in the problem.

A particularly important question concerns the relationship between Lepage
equivalents and null Lagrangians:\ that is, Lagrangians whose Euler-Lagrange
forms vanish. Clearly if a Lagrangian has a closed Lepage equivalent then it
will be null; and when $m = 1$ then the unique Lepage equivalent of a null
Lagrangian is closed. But when $m > 1$ then a choice of Lepage equivalent
would be needed, and it is not immediately obvious how this choice should be
made.

An answer to this question for first-order Lagrangians was found by
Krupka~\cite{Krup}, and also subsequently by Betounes~\cite{Bet}. In
coordinates $x^i$ on $M$ and fibred coordinates $(x^i, u^a)$ on $E$, the Lepage
equivalent
\[
\theta = L\,\omega
+ \sum_{r=1}^{\min\{m,n\}} \frac{1}{(r!)^2}
\frac{\partial^r L}{\partial u^{a_1}_{i_1} \ldots
\partial u^{a_r}_{i_r}}
\theta^{a_1} \wedge \ldots \wedge \theta^{a_r}
\wedge \omega_{i_1 \cdots i_r},
\]
of a Lagrangian $L \omega$ (where $\omega = dx^1 \wedge \ldots \wedge dx^m$ and
$\omega_{i_1 \cdots i_r} = i_{\partial / \partial x^{i_r}} \omega_{i_1 \cdots i_{r-1}}$,
and where $\theta^a = du^a - u^a_i dx^i$) is closed precisely when
$L \, \omega$ is null. We say that a Lepage equivalent constructed according to
this formula satisfies the \emph{closure property}. 
This property is significant in the context of symmetries 
because, as was pointed out in~\cite{Bet}, it allows us to conclude that any symmetry
of an arbitrary first-order Lagrangian is also a symmetry of its corresponding
form $\theta$, whereas this need not be the case for other Lepage equivalents.
An illustration of this phenomenon is given by Example~2 of that paper, where 
symmetries of the electromagnetic Lagrangian on $\R^4$ are considered. 
It is therefore natural to ask whether it is possible to find Lepage equivalents
with a similar property for higher-order Lagrangians. To date, though, no such
formula has been found for second-order (or higher-order) Lagrangians, and even the 
existence of Lepage equivalents having this additional property is unclear.

In this paper we look at homogeneous problems, where the idea of a Lepage equivalent 
is not directly appropriate. These homogeneous problems are defined on a manifold $E$ 
without any given fibration over a space of independent variables, where the solution to
the variational problem is a submanifold with an orientation but without any
preferred parametrization. Instead of using jet bundles for these problems, the
Lagrangian is defined instead on the bundle of $k$-th order $m$-frames $\Fm^k E$
in the manifold~\cite{CS1} (this is also called the bundle of regular $k$-th
order $m$-velocities). The Lagrangian is a function $L$ rather than an
$m$-form, and is required to satisfy a certain homogeneity condition.

Homogeneous variational problems arise in geometry more directly than in physics:\
for instance Finsler geometry considers the homogeneous problem with $m=1$ and $k=1$,
and the study of minimal surfaces involves a homogeneous problem with $m=2$. There
is, nevertheless, a straightforward relationship between the two types of problem, and the
paradigm of this is the `homogenization trick' of classical mechanics. Given a 
time-dependent Lagrangian 1-form $L(t, q^a, \dot{q}^a) dt$ defined on the jet bundle $J^1\pi$
where $\pi : \R \times M \to \R$, the corresponding homogeneous Lagrangian function is
\[
\widetilde{L}(t, \dot{t}, q^a, \dot{q}^a) = \dot{t} L \left( t, q^a, \dot{t}^{-1}\dot{q}^a \right) 
\]
defined on a suitable open subset of $T(\R \times M)$.

A similar relationship between the two types of problem holds in the general case.
Factoring the bundle of $m$-frames by the vector fields used to specify the homogeneity
condition gives rise to the bundle $J^k_+(E, m)$ of $k$-th order oriented
contact elements of dimension $m$; a Lagrangian $m$-form $\lambda$ on this bundle
gives rise to a homogeneous function $L$ on the frame bundle. If a fibration
$\pi : E \to M$ is given then there is an inclusion $J^k\pi \subset J^k_+(E, m)$,
and a Lagrangian form on $J^k\pi$ gives rise to a homogeneous function $L$ on 
an open subset of the frame bundle. More details of the relationship between
the two types of problem in the general case may be found in~\cite{CS1}, and given 
this relationship it is of some interest to search for $m$-forms related to
homogeneous Lagrangians having the closure property described above.

It was shown in~\cite{CS1} that, given any Lagrangian function $L$ on $\Fm^k E$,
it is possible to construct an $m$-form on $\Fm^{2k-1} E$ called the {\em
Hilbert-Carath\'{e}odory form}\/ having the same extremals as $L$ and giving
rise to a suitable Euler-Lagrange form. The Hilbert-Carath\'{e}odory form is
projectable to the bundle of contact elements when $m = 1$ or $k \leq 2$.
This form does not, in general, have the closure property. 
It was subsequently shown in~\cite{CS2} that, in the case of a first-order Lagrangian
function, there is another $m$-form on $\Fm^1 E$ with the property that it is
closed precisely when the Lagrangian is null. This second $m$-form is
projectable to the bundle of first-order contact elements, and if there is a 
fibration of $E$ over some $m$-dimensional manifold then the restriction to the
corresponding jet bundle takes the coordinate form shown above. We call this
second $m$-form the \emph{fundamental form} of the Lagrangian.

The present paper is a report on the second stage of a project to generalize the latter
construction to Lagrangians of arbitrary order:\ we have a candidate for the fundamental form, and our task is to show that it satisfies the closure property. The construction involves a sequence $\Theta_0, \Theta_1, \ldots, \Theta_m$ of vector-valued forms, where $\Theta_r$ is an $\bigwedge^{m-r}\R^{m*}$-valued $r$-form; $\Theta_m$ is thus a scalar $m$-form, and is our candidate. Our strategy for proving the closure property is to show that that $\Theta_r$ may be obtained from $\Theta_{r+1}$ by contraction with total derivatives. The proof of this for $r=0$ is straightforward; the proof for $r=1$ forms the significant content of this paper, and involves a substantial level of complexity. Combining this with a general result about such contractions shows that, in the case of a homogeneous Lagrangian of arbitrary order in two independent variables, the scalar $2$-form $\Theta_2$ satisfies the closure property. A previous note~\cite{Sau2Order} gave a version of this construction for second-order Lagrangians in two variables.

In Sections~2 and~3 of this paper we therefore recall the properties of homogeneous variational problems in the context of a bicomplex of vector-valued forms, and in Section~4 we collect together some preliminary results. The main theorem of the paper is in Section~5, and we offer a discussion of some consequences of this result in Section~6. We also suggest that, subject to overcoming the computational difficulties, it should be possible to extend the present results to Lagrangians in arbitrarily many independent variables.

\section{Homogeneous variational problems}

We consider a smooth manifold $E$ of dimension $n$, and its bundles $\tau_k : \Fm^k E \to E$ of $k$-th order $m$-frames. Important objects defined intrinsically on these bundles are the total derivatives and the vertical endomorphisms. The former are vector fields $\mathrm{T}_j$ along the map $\tau_{k+1,k} : \Fm^{k+1} E \to \Fm^k E$, and are described in coordinates as
\[
\mathrm{T}_j = \sum_{|I|=0}^k u^\alpha_{I+1_j} \vf{u^\alpha_I} \, ,
\]
and the latter are type $(1,1)$ tensor fields $S^i$ on $\Fm^{k+1} E$ described in coordinates as
\[
S^i = \sum_{|I|=0}^k \vf{u^\alpha_{I+1_i}} \otimes du^\alpha_I \, .
\]
Here and subsequently we use local coordinates $(u^\alpha)$ on $E$ and the corresponding jet coordinates $(u^\alpha_I)$ on $\Fm^k E$, where $I$ is a multi-index. Intrinsic definitions of the operators $\mathrm{T}_j$ and $S^i$ may be found in~\cite{CS1, SauHom}.

We also need to use the fundamental vector fields $\Delta^I_j$ defined by
\[
\Delta^I_j = S^I(\mathrm{T}_j)
\]
where the tensor fields $S^i$ and $S^j$ commute, so that $S^I$ may be defined by iteration; these vector fields are well-defined on the manifold $\Fm^{k+1} E$ (rather than along the map $\tau_{k+1,k}$).

We shall let $i_j$ denote the action of the vector field $\mathrm{T}_j$ on differential forms by contraction, and let $d_j$ denote the action as a Lie derivative; $i^I_j$ and $d^I_j$ will, similarly, denote the actions of the vector field $\Delta^I_j$. We shall use the symbol $S^i$ to denote the contraction of the tensor field with a form, as well as denoting the tensor field itself. The symbol $S^I$ will denote the iterated contraction, and we shall write $\St^I$ to denote the (single) contraction of the composite tensor field with the form; in the case of the action on a 1-form these are the same. 

As was demonstrated in~\cite{CS1}, a Lagrangian function $L$ on $\Fm^k E$ whose extremals have no preferred parametrization must be \emph{homogeneous}, in that it must satisfy the properties
\[
d^i_j L = \delta^i_j L \, , \qquad d^I_j L = 0 \quad \mbox{for } |I| \geq 2 \, .
\]
Associated with such a Lagrangian are its $m$ Hilbert forms. These are the 1-forms $\vartheta^i$ on $\Fm^{2k-1} E$ defined by
\[
\vartheta^i = \sum_{|I|=0}^k \frac{(-1)^{|I|}}{I! (|I|+1)} d_I S^{I+1_i} dL = P^i dL
\]
and are generalisations of the Hilbert form used in Finsler geometry. We shall need to use several properties of the Hilbert forms, and we record these below.

\textbf{Lemma 2.1}\\
\textit{The Hilbert forms $\vartheta^i$ have the following properties:
\begin{align*}
i_k \vartheta^i & = \delta^k_i L \, , \\
i^K_k \vartheta^i & = 0 \quad \text{when $|K| > 0$} \quad \text{and} \\
d^I_i \vartheta^j & = - \sum_M \frac{(-1)^{|M|}}{M!} C_{I,M,i,j} d_M S^{I + M - 1_i +1_j} dL
\end{align*}
where the coefficient $C_{I,M,i,j}$ is given by
\begin{align*}
C_{I,M,i,j} & = M(i) \left(
\frac{|I|! |M|! + (-1)^{|I|} (|I| + |M| - 1)!}{(|I| + |M| + 1)!} \right) \\
& \qquad - \, I(i) \left(
\frac{(|I| - 1)!(|M| + 1)! - (-1)^{|I|}(|I| + |M| - 1)!}{(|I| + |M| + 1)!}
\right) \\
& \qquad + \, \delta^j_i \left(
\frac{|I|! |M|! - (-1)^{|I|} (|I| + |M|)!}{(|I| + |M| + 1)!} \right) \, .
\end{align*}
}

\textbf{Proof}\\
These results are all derived in~\cite{CS1}: they are Lemma~5.5, Proposition~6.1 and the calculation immediately preceding Theorem~6.3.
\qed

The Hilbert forms are used to construct the Euler-Lagrange form
\[
\varepsilon = dL - d_i \vartheta^i
\]
on $\Fm^{2k} E$. In coordinates
\[
\varepsilon = \sum_{|I|=0}^k (-1)^{|I|} d_I \left( \pd{L}{u^\alpha_I} \right) du^\alpha
\]
incorporating the Euler-Lagrange equations for the variational problem defined by $L$. More details of this construction may be found in~\cite{CS1}.

The fact that we have a family of Hilbert forms for multiple-integral problems suggests that it might be advantageous to consider them as the components of a vector-valued form. Spaces of suitable vector-valued forms were introduced in~\cite{SauHom}; these are the spaces
\[
\Omega^{r,s}_k = \Omega^r (\Fm^k E) \otimes {\textstyle \bigwedge^s} \R^{m*}
\]
of $r$-forms on $\Fm^k E$ taking their values in the vector space of alternating $s$-linear forms on $\R^m$. We denote the standard basis of $\R^{m*}$ suggestively by $(dt^i)$, so that an element $\Phi \in \Omega^{r,s}_k$ would be written in components as $\phi_{i_1 \cdots i_s} \otimes dt^{i_1} \wedge \cdots \wedge dt^{i_s}$, where the $\phi_{i_1 \cdots i_s}$ are scalar $r$-forms completely skew-symmetric in their indices.

A significant feature of these spaces is that they may be used to form a family of \emph{variational bicomplexes}, rather like the variational bicomplex of scalar forms on jet bundles. The mappings between the spaces are the ordinary de~Rham differential $d : \Omega^{r,s}_k \to \Omega^{r+1,s}_k$ acting on the individual components of the form, and the \emph{total derivative operator} $\dT : \Omega^{r,s}_k \to \Omega^{r,s+1}_{k+1}$ defined by
\[
\dT \left( \phi_{i_1 \cdots i_s} \otimes dt^{i_1} \wedge \cdots \wedge dt^{i_s} \right)
= \left( d_j \phi_{i_1 \cdots i_s} \right) 
\otimes dt^j \wedge dt^{i_1} \wedge \cdots \wedge dt^{i_s} \, .
\]
The initial part of such a bicomplex is shown in the diagram below,
\begin{figure}
\setlength{\unitlength}{0.9pt}
\begin{center}
\begin{picture}(450,490)(0,-170)
\multiput(160,300)(80,0){3}{\vector(0,-1){40}}
\multiput(160,220)(80,0){3}{\vector(0,-1){40}}
\multiput(160,120)(80,0){3}{\makebox(0,0){$\vdots$}}
\multiput(160,60)(80,0){3}{\vector(0,-1){40}}
\multiput(160,-20)(80,0){3}{\vector(0,-1){40}}
\multiput(160,-100)(80,0){3}{\vector(0,-1){40}}
\multiput(105,240)(0,-80){2}{\vector(1,0){30}}
\put(110,80){\vector(1,0){25}}
\multiput(105,0)(0,-80){2}{\vector(1,0){30}}
\multiput(180,240)(0,-80){5}{\vector(1,0){35}}
\put(180,-20){\vector(1,-1){40}}
\multiput(260,240)(0,-80){5}{\vector(1,0){35}}
\multiput(360,240)(0,-80){5}{\makebox(0,0){$\cdots$}}
\multiput(80,240)(0,-80){5}{\makebox(0,0){$0$}}
\put(160,320){\makebox(0,0){$0$}}
\multiput(240,320)(80,0){2}{\makebox(0,0){$0$}}
\put(160,240){\makebox(0,0){$\Obar^{0,0}_k$}}
\put(240,240){\makebox(0,0){$\Omega^{1,0}_k$}}
\put(320,240){\makebox(0,0){$\Omega^{2,0}_k$}}
\put(160,160){\makebox(0,0){$\Obar^{0,1}_{k+1}$}}
\put(240,160){\makebox(0,0){$\Omega^{1,1}_{k+1}$}}
\put(320,160){\makebox(0,0){$\Omega^{2,1}_{k+1}$}}
\put(160,80){\makebox(0,0){$\Obar^{0.m-1}_{k+m-1}$}}
\put(240,80){\makebox(0,0){$\Omega^{1,m-1}_{k+m-1}$}}
\put(320,80){\makebox(0,0){$\Omega^{2,m-1}_{k+m-1}$}}
\put(160,0){\makebox(0,0){$\Obar^{0,m}_{k+m}$}}
\put(240,0){\makebox(0,0){$\Omega^{1,m}_{k+m}$}}
\put(320,0){\makebox(0,0){$\Omega^{2,m}_{k+m}$}}
\put(160,-80){\makebox(0,0){$\overline{Xi}^0_{k+m}$}}
\put(240,-80){\makebox(0,0){$\Xi^1_{k+m}$}}
\put(320,-80){\makebox(0,0){$\Xi^2_{k+m}$}}
\multiput(160,-160)(80,0){3}{\makebox(0,0){$0$}}
\multiput(200,10)(0,80){4}{\makebox(0,0)[b]{$\scriptstyle d$}}
\multiput(280,10)(0,80){4}{\makebox(0,0)[b]{$\scriptstyle d$}}
\multiput(200,-70)(80,0){2}{\makebox(0,0)[b]{$\scriptstyle \delta$}}
\multiput(155,200)(80,0){3}{\makebox(0,0)[r]{$\scriptstyle\dT$}}
\multiput(155,40)(80,0){3}{\makebox(0,0)[r]{$\scriptstyle\dT$}}
\put(155,-40){\makebox(0,0)[r]{$\scriptstyle p_0$}}
\put(235,-40){\makebox(0,0)[r]{$\scriptstyle p_1$}}
\put(315,-40){\makebox(0,0)[r]{$\scriptstyle p_2$}}
\end{picture}
\end{center}
\setlength{\unitlength}{1pt}
\caption{The homogeneous variational bicomplex}
\end{figure}
where we have written $\Xi^r_k$ for the quotient $\Omega^{r,m}_k / \dT(\Omega^{r,m-1}_{k-1})$, and used an overline for the spaces of vector-valued functions in the first column to denote quotients by the constant functions. Of course the forms in row $s$ of a full bicomplex must necessarily be defined on $\Fm^k E$ where $k \geq s$, but there are also partial bicomplexes which omit the top part of the diagram, starting in row $s$ with forms defined on $E$, and finishing in row $m$ with forms defined on $\Fm^{m-s} E$.

The rows of each bicomplex are of course locally exact. The columns of each individual bicomplex are not exact, even locally; but if we work modulo pull-backs then each column apart from the first is in fact globally exact:\ that is, if the vector-valued form $\Phi \in \Omega^{r,s}_k$ satisfies $\dT \Phi = 0$ then there is a form $\Psi \in \Omega^{r,s-1}_l$ with $l \geq k$ such that $\dT \Psi = \Phi$. By diagram chasing, therefore, the first column is also locally exact modulo pull-backs. The homotopy operator used to show global exactness in this sense is the map $P : \Omega^{r,s}_k \to \Omega^{r,s+1}_{(r+1)k-1}$ defined (see~\cite{SauHom}) by
\begin{align*}
P\Phi & = P^j_{(s)}(\phi_{i_1 \cdots i_s}) \otimes
\left\{ \vf{t^j} \hook \left( dt^{i_1} \wedge \ldots
\wedge dt^{i_s} \right) \right\} \\
& = s \, P^j_{(s)}(\phi_{j i_2 \cdots i_s})
\otimes dt^{i_2} \wedge \ldots \wedge dt^{i_s}
\end{align*}
where $P^j_{(s)}$ is the differential operator on scalar $r$-forms defined by
\[
P^j_{(s)} = \sum_{|J| = 0}^{rk-1}
\frac{(-1)^{|J|}(m - s)! |J|!}{r^{|J|+1} (m - s + |J| + 1)! J!}
d_J S^{J + 1_j} \, .
\]

The homogeneous variational problems described above fit comfortably within this framework, where we use the notation
\[
d^m t = dt^1 \wedge \cdots \wedge dt^m \, ,
\qquad d^{m-1}t_j = \vf{t^j} \hook d^m t \, .
\]
If we let the Hilbert forms $\vartheta^j$ be the components of a vector-valued 1-form $\Theta_1 \in \Omega^{1,m-1}_{2k-1}$ so that $\Theta_1 = \vartheta^j \otimes d^{m-1}t_j$, and if we let $\Theta_0 \in \Obar^{0,m}_k$ be the equivalence class of the vector-valued function $L \, d^m t$, then the formula for the Hilbert forms may be written simply as $\Theta_1 = Pd\Theta_0$. If, similarly, we let $\ce_0$ denote the vector-valued 1-form $\varepsilon \otimes d^m t$ then the Euler-Lagrange formula is simply $\ce_0 = d\Theta_0 - \dT\Theta_1$.

\section{The fundamental form of a homogeneous Lagrangian}

The simplicity of the formula $\Theta_1 = Pd\Theta_0$ for the vector of Hilbert forms suggests that we might wish to consider higher powers of the operator $Pd$ acting on the Lagrangian. Put
\[
\left.
\begin{array}{rcl}
\Theta_q & = & (Pd)^q \Theta_0 \in \O^{q,m-q} \\[1ex]
\ce_q & = & (Pd)^q \ce_0 \in \O^{q+1,m-q}
\end{array}
\right\}
\qquad 0 \leq q \leq m \, ,
\]
where we have omitted explicit mention of the order of the manifold on which these forms are defined as its complicated expression tends to obscure the overall message.

\textbf{Lemma 3.1}
\[
\ce_q = (-1)^q (d\Theta_q - \dT \Theta_{q+1}) \qquad 0 \leq q \leq m-1 \, .
\]
\textbf{Proof}
If this relationship holds for some given value of $q$ then
\begin{multline*}
\ce_{q+1} = Pd\ce_q = (-1)^{q+1} Pd\dT\Theta_{q+1} 
= (-1)^{q+1} P\dT d\Theta_{q+1} \\
= (-1)^{q+1} (d\Theta_{q+1} - \dT Pd\Theta_{q+1}) 
= (-1)^{q+1} (d\Theta_{q+1} - \dT\Theta_{q+2}) \, ,
\end{multline*}
and the relationship is certainly true for $q=0$. \qed

Our main interest will be in $\Theta_m \in \Omega^{m,0}$ as this takes its values in the 1-dimensional vector space $\bigwedge^0 \R^{m*}$ and may be identified with a scalar $m$-form; we shall call this the \emph{fundamental form} of the Lagrangian $L$. If we consider the case of a first-order Lagrangian, we can use the formula for $P$ to give an explicit description of the fundamental form:\ in this case each $\Theta_q$ is also first order, and we see easily that
\[
\Theta_m = \frac{1}{m!} (S^1 d) \ldots (S^m d) L \, ;
\]
this is just the $m$-form described in~\cite{CS2} and shown there to satisfy the closure property. It is also shown that this form, defined on $\Fm^1 E$, is projectable to the bundle of contact elements, and that if $L$ is derived by homogenisation from a Lagrangian on a jet bundle then the projection of $\Theta_m$ is just the fundamental Lepage equivalent found by Krupka and Betounes.

It is now natural to ask whether a similar property holds for higher-order Lagrangians. The construction of $\Theta_m$ may be carried out for a Lagrangian of arbitrary order, and so we are led to the following conjecture.

\textbf{Conjecture}
\emph{Let $L$ be a homogeneous Lagrangian defined on $\Fm^k E$; then $L$ is null if, and only if, $d\Theta_m = 0$.}

In one direction the proof is straightforward, and has nothing to do with homogeneity. Suppose that $\ce_0 = 0$, so that $d\Theta_0 = \dT \Theta_1$. Then, recursively,
\[
d\Theta_q = \dT\Theta_{q+1} \qquad 0 \leq q \leq m-1 \, ;
\]
for if this relationship holds then
\[
P\dT d\Theta_{q+1} = Pd\dT\Theta_{q+1} = Pd^2 \Theta_q = 0
\]
so that
\[
d\Theta_{q+1} = \dT Pd\Theta_{q+1} = \dT\Theta_{q+2}
\]
using the homotopy formula and the definition of $\Theta_{q+2}$. In particular $d\Theta_{m-1} = \dT\Theta_m$; so finally, therefore,
\[
d\Theta_m = P\dT d\Theta_m = Pd\dT\Theta_m = Pd^2 \Theta_{m-1} = 0 \, .
\]

The converse, that the closure of $\Theta_m$ implies the nullity of $L$, is much harder, and homogeneity is essential:\ for instance, in the single-integral case, take a non-zero Lagrangian depending on only the position coordinates. Such a Lagrangian is certainly not null, but $\Theta_1 = 0$ so that $\Theta_1$ is certainly closed.

As a first step towards a proof, suppose that the \emph{recovery formula} holds:\ that is, that
\[
\Theta_q = \frac{1}{m-q} \, \iT \Theta_{q+1} \, , \qquad 0 \leq q \leq m-1 \, ,
\]
where $\iT$ denotes contraction with the total derivative operator,
\[
\iT \left( \phi_{i_1 \cdots i_s} \otimes dt^{i_1} \wedge \cdots \wedge dt^{i_s} \right)
= \left( i_j \phi_{i_1 \cdots i_s} \right) 
\otimes dt^j \wedge dt^{i_1} \wedge \cdots \wedge dt^{i_s} \, ;
\]
note that $\dT = d\iT + \iT d$ and that $\iT \dT + \dT \iT = 0$.

\textbf{Lemma 3.2}\\
\textit{If the recovery formula holds then}
\[
\ce_q = \frac{1}{m-q} \, \iT\ce_{q+1} \, .
\]
\textbf{Proof}
\begin{align*}
\iT\ce_{q+1} & = (-1)^{q+1} (\iT d\Theta_{q+1} - \iT \dT \Theta_{q+2}) \\
& = (-1)^{q+1} (\dT\Theta_{q+1} - d\iT\Theta_{q+1} + \dT \iT \Theta_{q+2}) \\
& = (-1)^{q+1} (\dT\Theta_{q+1} - (m-q)d\Theta_q + (m-q-1)\dT \Theta_{q+1}) \\
& = (-1)^q (m-q) (d\Theta_q - \dT\Theta_{q+1} ) \\
& = (m-q) \ce_q \, .
\end{align*}
\qed

\textbf{Lemma 3.3}\\
\textit{If the recovery formula holds and $d\Theta_m = 0$ then the Lagrangian is null.}

\textbf{Proof}
\begin{align*}
\ce_{m-1} & = (-1)^{m-1} (d\Theta_{m-1} - \dT\Theta_m) \\
& = (-1)^{m-1} (d\iT\Theta_m - (d\iT\Theta_m + \iT d\Theta_m) ) \\
& = (-1)^m \iT d\Theta_m \\
& = 0 \, .
\end{align*}
It now follows from Lemma~3.2 that
\[
\ce_0 = \frac{1}{m!} \, \ce_{m-1} = 0
\]
so that the Lagrangian is null. \qed

We are therefore led to the question of whether the recovery formula holds for a general homogeneous Lagrangian. The rest of this paper is devoted to proving that the first two steps hold, so that
\[
\Theta_0 = \frac{1}{m} \, \iT\Theta_1 \, , \qquad \Theta_1 = \frac{1}{m-1} \, \iT\Theta_2 \, .
\]

\section{Some preliminary results}

In order to achieve our objective, we need to consider the scalar components of the vector-valued forms $\iT\Theta_1$ and $\iT\Theta_2$. It is convenient to introduce the operator $\DD_p$ on scalar forms, defined by
\[
\DD_p = \sum_{|I|=p} \frac{1}{I!} d_I S^I \, .
\]
For some calculations we shall need to use the equivalent expression for $\DD_p$ using a list of ordinary indices rather than a single multi-index, and this is
\[
\DD_p = \frac{1}{p!} d_{i_1 \cdots i_p} S^{i_1 \cdots i_p} \, ;
\]
the sum over the ordinary indices is understood. The conversion between the two types of notation involves, for a given multi-index $I$, the quantity $|I|!/I!$ known as its \emph{weight}; this quantity is the ratio between symmetrized and non-symmetrized index expressions.

Write
\[
\Theta_1 = \vartheta^j \otimes d^{m-1}t_j \, , \qquad
\Theta_2 = \vartheta^{ij} \otimes d^{m-2}t_{ij} 
\]
where
\[
\vartheta^j = P_{(1)}^j dL \, , \qquad \vartheta^{ij} = P_{(2)}^j \vartheta^i - P_{(2)}^i \vartheta^j \, ;
\]
we then have
\begin{align*}
P_{(1)}^j & = \sum_{p=0}^k \frac{(-1)^p}{p+1} \DD_p S^j \, , \\
P_{(2)}^j & = \sum_{p=0}^{2k-1} \frac{(-1)^p \, p!}{2^{p+1} (p+2)!} \DD_p S^j \, .
\end{align*}
Evaluating the contraction with a total derivative thus involves moving $i_j$ to the right of $\DD_p S^j$, so that we can use the homogeneity properties of $L$; we therefore need to consider commutators.

\textbf{Lemma 4.1} 

\textit{The commutators of total derivative operators and vertical endomorphisms are given by the following formul\ae:}
\begin{align*}
[i^I_i, d_j] & = I(j) i^{I - 1_j}_i &
[d^I_i, \St^J] & = -J(i) \, \St^{J+I-1_i} \\{}
[i^I_i, \St^J] & = i^{I+J}_i &
[d^I_i, d_j] & = I(j) d^{I - 1_j}_i \, .
\end{align*}
Proofs of these formul\ae\ may be found in~\cite{CS1}, or are easy consequences of the results there. \qed

\textbf{Lemma 4.2}
\[
S^J \DD_p = \sum_{q=0}^{|J|} \frac{|J|!}{q! \, (|J|-q)!} \DD_{p-q} S^J
\]
where we adopt the convention that $\DD_{p-q} = 0$ when $q > p$.

\textbf{Proof}
We use induction on the length of the multi-index $J$, and also~\cite[Lemma~2.1]{CS1} which in the present notation reads
\[
S^j \DD_p = (\DD_p + \DD_{p-1})S^j \, .
\]
So suppose the proposed formula is true for every multi-index of length $r$, and that $|K|=r+1$. Put $K=J+1_j$ so that $|J|=r$, and then
\begin{align*}
S^K \DD_p & = S^j S^J \DD_p \\
& = S^j \sum_{q=0}^{|J|} \frac{|J|!}{q! \, (|J|-q)!} \DD_{p-q} S^J \\
& = \sum_{q=0}^{|J|} \frac{|J|!}{q! \, (|J|-q)!} (\DD_{p-q} + \DD_{p-q-1}) S^j S^J \\
& = \sum_{q=0}^{|J|} \frac{|J|!}{q! \, (|J|-q)!} \DD_{p-q} S^K 
+ \sum_{q=0}^{|J|} \frac{|J|!}{q! \, (|J|-q)!} \DD_{p-q-1} S^K \\
& = \sum_{q=0}^{|J|} \frac{|J|!}{q! \, (|J|-q)!} \DD_{p-q} S^K 
+ \sum_{q=1}^{|J|+1} \frac{|J|!}{(q-1)! \, (|J|-q+1)!} \DD_{p-q} S^K \\
& = \sum_{q=0}^{|J|+1} \frac{|J|!((|J|-q+1)+q)}{q! \, (|J|-q+1)!} \DD_{p-q} S^K \\
& = \sum_{q=0}^{|K|} \frac{|K|!}{q! \, (|K|-q)!}  \DD_{p-q} S^K
\end{align*} 
as required; of course when $|J|=1$ this is just~\cite[Lemma~2.1]{CS1}. \qed

For convenience we shall put
\[
G_{|J|,q} = \frac{|J|!}{q! \, (|J|-q)!} \, ,
\]
and if we don't wish to use the convention regarding $\DD_{p-q}$ when $q > p$ then we simply write the sum as
\[
S^J \DD_p = \sum_{q=0}^{\min\{|J|,p\}} G_{|J|,q} \DD_{p-q} S^J \, .
\]

\textbf{Lemma 4.3}
\[
i_k \DD_p = \sum_{|K|=0}^p \sum_{|J|=p-|K|}
\frac{1}{J! \, K!} d_{J+K} S^J i^K_k \, .
\]
\textbf{Proof}

We carry out the proof using ordinary indices, and we claim that
\[
i_k \DD_p = d_{i_1} \cdots d_{i_p} 
\sum_{q=0}^p \frac{1}{q! \, (p-q)!} S^{i_1} \cdots S^{i_q} i^{i_{q+1} \cdots i_p}_k
\]
(summed, of course, over $i_1, \ldots, i_p$). For suppose that this is true for some value of $p$; then
\begin{align*}
i_k \DD_{p+1}
& = \frac{1}{p+1} i_k d_j \DD_p S^j 
= \frac{1}{p+1} d_j i_k \DD_p S^j \\
& = \frac{1}{p+1} d_j d_{i_1} \cdots d_{i_p} 
\sum_{q=0}^p \frac{1}{q! \, (p-q)!} S^{i_1} \cdots S^{i_q} i^{i_{q+1} \cdots i_p}_k S^j \\
& = \frac{1}{p+1} d_j d_{i_1} \cdots d_{i_p} 
\sum_{q=0}^p \frac{1}{q! \, (p-q)!} S^{i_1} \cdots S^{i_q} 
(S^j i^{i_{q+1} \cdots i_p}_k + i^{i_{q+1} \cdots i_p j}_k) \\
& = \frac{1}{p+1} d_{i_1} \cdots d_{i_{q+1}} \cdots d_{i_{p+1}} 
\sum_{q=0}^p \frac{1}{q! \, (p-q)!} S^{i_1} \cdots S^{i_q} 
S^{i_{q+1}} i^{i_{q+2} \cdots i_{p+1}}_k \\
& \qquad + \, \frac{1}{p+1} d_{i_1} \cdots d_{i_p} d_{i_{p+1}}
\sum_{q=0}^p \frac{1}{q! \, (p-q)!} S^{i_1} \cdots S^{i_q} 
i^{i_{q+1} \cdots i_p i_{p+1}}_k 
\end{align*}
where in the first sum we have relabelled the indices $i_{q+1}, \ldots, i_p$ as $i_{q+2}, \ldots, i_{p+1}$ and then relabelled the index $j$ as $i_{q+1}$, and in the second sum we have just relabelled the index $j$ as $i_{p+1}$. But now, in the first sum, replace $q$ by $q-1$ to give
\begin{align*}
i_k \DD_{p+1}
& = \frac{1}{p+1} d_{i_1} \cdots d_{i_{p+1}} 
\sum_{q=1}^{p+1} \frac{1}{(q-1)! \, (p+1-q)!} S^{i_1} \cdots S^{i_q} 
i^{i_{q+1} \cdots i_{p+1}}_k \\
& \qquad + \, \frac{1}{p+1} d_{i_1} \cdots d_{i_{p+1}}
\sum_{q=0}^p \frac{1}{q! \, (p-q)!} S^{i_1} \cdots S^{i_q} 
i^{i_{q+1} \cdots i_{p+1}}_k \\
& = \frac{1}{p+1} d_{i_1} \cdots d_{i_{p+1}} 
\sum_{q=0}^{p+1} \frac{q + (p+1-q)}{q! \, (p+1-q)!} S^{i_1} \cdots S^{i_q} 
i^{i_{q+1} \cdots i_{p+1}}_k \\
& = d_{i_1} \cdots d_{i_{p+1}} 
\sum_{q=0}^{p+1} \frac{1}{q! \, (p+1-q)!} S^{i_1} \cdots S^{i_q} 
i^{i_{q+1} \cdots i_{p+1}}_k 
\end{align*}
as required.

We now need to reinstate the multi-indices, and so finally we have
\begin{align*}
i_k \DD_p & = d_{i_1} \cdots d_{i_p} 
\sum_{q=0}^p \frac{1}{q! \, (p-q)!} S^{i_1} \cdots S^{i_q} i^{i_{q+1} \cdots i_p}_k \\
& = \sum_{q=0}^p \frac{1}{q! \, (p-q)!} \sum_{|K|=q} \sum_{|J|=p-q}
\frac{|J|! \, |K|!}{J! \, K!} d_{J+K} S^J i^K_k \\
& = \sum_{|K|=0}^p \sum_{|J|=p-|K|}
\frac{1}{J! \, K!} d_{J+K} S^J i^K_k \, .
\end{align*}
\qed

There is one further formula we need; this is a straightforward property of the Hilbert forms $\vartheta^i$.

\textbf{Lemma 4.4}
\[
S^i \vartheta^j = S^j \vartheta^i
\]

\textbf{Proof}\\
Using~\cite[Lemma~2.1]{CS1} in the form given above we have
\begin{align*}
S^i \vartheta^j
& = S^i \sum_{p=0}^k \frac{(-1)^p}{p+1} \DD_p S^j dL \\
& = \sum_{p=0}^k \frac{(-1)^p}{p+1} (\DD_p + \DD_{p-1}) S^i S^j dL \\
& = S^j \vartheta^i \, . 
\end{align*}
\qed

\section{The recovery formula}

We now give, as advertised, the proofs of steps 1 and 2 of the recovery formula. The proof of step 1 is, in fact, comparatively straightforward.

\textbf{Theorem 5.1}
\[
\Theta_0 = \frac{1}{m} \, \iT\Theta_1 \, .
\]
\textbf{Proof}

We have
\begin{align*}
\iT \Theta_1 & = i_i P^j dL \otimes dt^i \wedge d^{m-1}t_j \\
& = i_i P^j dL \otimes \delta^i_j d^m t \\
& = i_j P^j dL \otimes d^m t \, ,
\end{align*}
and
\[
i_j P^j dL =  \sum_{p=0}^k
\frac{(-1)^p}{p+1} i_j \DD_p S^j dL \, ;
\]
but from Lemmas 4.1 and 4.2 we have 
\begin{align*}
i_j \DD_p S^j dL & = \sum_{|K|=0}^p \sum_{|J|=p-|K|}
\frac{1}{J! \, K!} d_{J+K} S^J i^K_k S^j dL \\
& = \sum_{|K|=0}^p \sum_{|J|=p-|K|}
\frac{1}{J! \, K!} d_{J+K} S^J (S^j i^K_j + i^{K+1_j}_j) dL \, .
\end{align*}
Now $S^J i^K_j dL = 0$ because $i^K_j dL$ is a 0-form, and $i^{K+1_j}_j dL = 0$ unless $K=0$ by homogeneity. But then $S^J i^j_j dL = 0$ unless $J=0$, so we have $i_j \DD_p S^j dL = 0$ when $p>0$. We conclude that
\[
\iT \Theta_1 = i_j P^j dL \otimes d^m t = i^j_j dL \otimes d^m t 
= m L \otimes d^m t = \Theta_0 \, .
\]
\qed

\textbf{Corollary 5.2}\\
\textit{The fundamental form $\Theta_1$ of a single-integral homogeneous Lagrangian $L$ satisfies the closure condition.} \qed

The second step is, however, considerably more complicated.

\textbf{Theorem 5.3}
\[
\Theta_1 = \frac{1}{m-1} \, \iT\Theta_2 \, .
\]
\textbf{Proof}
We must show that
\[
i_j \left( P^j_{(2)} d\vartheta^i - P^i_{(2)} \vartheta^j \right) = (m-1) \vartheta^i
\quad\text{where}\quad
P^j_{(2)} = \sum_p \frac{(-1)^p \, p!}{2^{p+1} (p+2)!} \DD_p S^j \, ,
\]
in other words that
\[
i_j \left( \sum_p \lambda_p \DD_p S^j d\vartheta^i 
- \sum_p \lambda_p \DD_p S^i \vartheta^j \right) = (m-1) \vartheta^i 
\quad\text{where}\quad
\lambda_p = \frac{(-1)^p \, p!}{2^{p+1} (p+2)!} \, .
\]
Now
\begin{align*}
i_k \sum_p \lambda_p \DD_p S^j d\vartheta^i
& = \sum_p \lambda_p i_k \DD_p S^j d\vartheta^i \\
& = \sum_p \lambda_p \sum_{|K|=0}^p \sum_{|J|=p-|K|}
\frac{1}{J! \, K!} d_{J+K} S^J i^K_k S^j d\vartheta^i \\
& = \sum_p \lambda_p \sum_{|K|=0}^p \sum_{|J|=p-|K|}
\frac{1}{J! \, K!} d_{J+K} (S^{J+1_j} i^K_k d\vartheta^i + S^J i^{K+1_j}_k d\vartheta^i) \\
& = \sum_p \lambda_p \sum_{|K|=0}^p \sum_{|J|=p-|K|}
\frac{1}{J! \, K!} d_{J+K} (S^{J+1_j} d^K_k \vartheta^i + S^J d^{K+1_j}_k \vartheta^i) \\
& \qquad - \, \sum_p \lambda_p \sum_{|J|=p}
\frac{1}{J!} d_J S^{J+1_j} d(\delta^i_k L)
\end{align*}
using Lemma~2.1; we consider the contributions of the two sums separately. The second sum is just
\[
- \delta^i_k \sum_p \lambda_p \DD_p S^j dL
\]
and so in the required skew combination gives
\begin{align}
\hspace{-8em}\makebox[0em]{$\displaystyle - \delta^i_j \sum_p \lambda_p \DD_p S^j dL
+ \delta^j_j \sum_p \lambda_p \DD_p S^i dL$} \nonumber \\
& = (m-1) \sum_p \lambda_p \DD_p S^i dL \, ; \label{part1}
\end{align}
it is the first sum which will require more careful attention. We split it into two parts. The first, in the required skew combination, is
\[
\sum_p \lambda_p \sum_{|K|=0}^p \sum_{|J|=p-|K|}
\frac{1}{J! \, K!} d_{J+K} (S^{J+1_j} d^K_j \vartheta^i - S^{J+1_i} d^K_j \vartheta^j) \, ;
\]
but $S^j d^K_k = d^K_k S^j + \delta^j_k \St^K$ so that
\begin{align*}
S^j d^K_j \vartheta^i - S^i d^K_j \vartheta^j
& = (d^K_j S^j + \delta^j_j \St^K) \vartheta^i - (d^K_k S^j + \delta^i_j \St^K) \vartheta^j \\
& = d^K_j (S^j \vartheta^i - S^i \vartheta^j) + (m-1) \St^K \vartheta^i \\
& = (m-1) \St^K \vartheta^i
\end{align*}
using Lemma~4.4, so this part becomes
\[
(m-1) \sum_p \lambda_p \sum_{|K|=0}^p \sum_{|J|=p-|K|}
\frac{1}{J! \, K!} d_{J+K} S^J \St^K \vartheta^i \, .
\]
In fact, as $\vartheta^i$ is a 1-form, we can replace $\St^K$ by $S^K$ to obtain
\begin{align*}
\lefteqn{\hspace{-2em}(m-1) \sum_p \lambda_p \sum_{|K|=0}^p \sum_{|J|=p-|K|}
\frac{1}{J! \, K!} d_{J+K} S^J S^K \vartheta^i} \\
& =  (m-1) \sum_p \lambda_p \sum_{q=0}^p 
\frac{1}{q! \, (p-q)!} d_{k_1 \cdots k_p} S^{k_1 \cdots k_p} \vartheta^i \\
& = (m-1) \sum_p 
\frac{2^p \lambda_p}{p!} d_{k_1 \cdots k_p} S^{k_1 \cdots k_p} \vartheta^i \\
& = (m-1) \sum_p 
\frac{(-1)^p}{2(p+2)!} d_{k_1 \cdots k_p} S^{k_1 \cdots k_p} \vartheta^i \\
& = (m-1) \sum_K 
\frac{(-1)^{|K|} |K|!}{2(|K|+2)! \, K!} d_K S^K \vartheta^i \, .
\end{align*}

We now, of course, need to replace $S^K \vartheta^i$ by its expression in terms of $dL$. We get
\[
(m-1) \sum_K 
\frac{(-1)^{|K|} |K|!}{2(|K|+2)! \, K!} d_K 
\sum_M \frac{(-1)^{|M|} |M|! \, |K|!}{(|M|+|K|+1)! \, M!} d_M S^{M+K+1_i} dL 
\]
and so we need to replace the double sum by a single one. As before, we transform the multi-indices into ordinary indices before attempting this. We have
\begin{align*}
\lefteqn{\hspace{-2em}(m-1) \sum_K \sum_M 
\frac{(-1)^{|K|+|M|} |K|!}{2(|K|+2)! \, (|M|+|K|+1)!} 
\frac{|K|! \, |M|!}{K! \, M!} d_{K+M} S^{K+M+1_i} dL} \\
& = (m-1) \sum_{p,q} \sum_{|K|=p} \sum_{|M|=q} 
\frac{(-1)^{p+q} p!}{2(p+2)! \, (p+q+1)!} 
\frac{|K|! \, |M|!}{K! \, M!} d_{K+M} S^{K+M+1_i} dL \\
& = (m-1) \sum_{p=0}^\infty \sum_{q=0}^\infty \frac{(-1)^{p+q} p!}{2(p+2)! \, (p+q+1)!}
d_{k_1 \cdots k_p k_{p+1} \cdots k_{p+q}} S^{k_1 \cdots k_p k_{p+1} \cdots k_{p+q} i} dL \\
& = (m-1) \sum_{r=0}^\infty \sum_{p=0}^r \frac{(-1)^r p!}{2(p+2)! \, (r+1)!}
d_{k_1 \cdots k_r} S^{k_1 \cdots k_r i} dL 
\end{align*}
on replacing $p+q$ by $r$; but
\[
\frac{p!}{(p+2)!} = \frac{1}{(p+1)(p+2)} = \frac{1}{p+1} - \frac{1}{p+2}
\]
so that
\[
\sum_{p=0}^r \frac{p!}{(p+2)!} = 1 - \frac{1}{r+2} = \frac{r+1}{r+2} \, .
\]
We are therefore left with
\[
(m-1) \sum_{r=0}^\infty \frac{(-1)^r (r+1)}{2 \, (r+2)!}
d_{k_1 \cdots k_r} S^{k_1 \cdots k_r i} dL \, , 
\]
and restoring the multi-indices gives
\[
(m-1) \sum_{|K|=0}^\infty \frac{(-1)^{|K|}}{2 \, (|K|+2) \, K!}
d_K S^{K+1_i} dL 
\]
which we may equivalently write as
\begin{equation}
(m-1) \sum_{p=0}^\infty \frac{(-1)^p}{2 \, (p+2)} \DD_p S^i dL \, . \label{part2}
\end{equation}

We now move on to the final part of the expression we are analysing. The third formula in Lemma~2.1 gives
\begin{align*}
\lefteqn{\hspace{-2em}d^{K+1_j}_j \vartheta^i - d^{K+1_i}_j \vartheta^j} \\ 
& = \sum_M \frac{(-1)^{|M|}}{M!} \left( C_{K+1_i,M,j,j} d_M S^{M + K + 1_i} dL 
- C_{K+1_j,M,j,i} d_M S^{M + K + 1_i} dL \right) \\
& = \sum_M \frac{(-1)^{|M|}}{M!} 
\left( C_{K+1_i,M,j,j} - C_{K+1_j,M,j,i} \right) d_M S^{M + K + 1_i} dL \, ,
\end{align*}
and from the expression for $C_{I,M,i,j}$ we obtain
\[
C_{K+1_i,M,j,j} - C_{K+1_j,M,j,i}
= (m-1) \left(
\frac{|K|! |M|! + (-1)^{|K|} (|K| + |M|)!}{(|K| + |M| + 1)!} \right) \, .
\]
We shall write $(m-1) F_{|K|,|M|}$ for this latter coefficient, so that
\[
d^{K+1_j}_j \vartheta^i - d^{K+1_i}_j \vartheta^j 
= (m-1) \sum_M \frac{(-1)^{|M|}}{M!} 
F_{|K|,|M|} d_M S^{K + M + 1_i} dL \, .
\]
Thus
\begin{align*}
\lefteqn{\hspace{-2em}\sum_p \lambda_p \sum_{|K|=0}^p \sum_{|J|=p-|K|}
\frac{1}{J! \, K!} d_{J+K} S^J (d^{K+1_j}_j \vartheta^i - d^{K+1_i}_j \vartheta^j)} \\
& = (m-1) \sum_p \lambda_p \sum_{|K|=0}^p \sum_{|J|=p-|K|} \sum_M 
\frac{(-1)^{|M|}}{J! \, K! \, M!} F_{|K|,|M|} 
d_{J+K} S^J d_M S^{K + M + 1_i} dL
\end{align*}
and we need to move $S^J$ to the right of $d_M$. We get
\begin{align*}
\lefteqn{\hspace{-1em}(m-1) \sum_p \lambda_p \sum_{|K|=0}^p \sum_{|J|=p-|K|} \sum_M 
\frac{(-1)^{|M|}}{J! \, K! \, M!} F_{|K|,|M|} 
d_{J+K} S^J d_M S^{K + M + 1_i} dL} \\
& = (m-1) \sum_p \lambda_p \sum_{|K|=0}^p \sum_{|J|=p-|K|} \sum_r 
\frac{(-1)^r}{J! \, K!} F_{|K|,r} 
d_{J+K} S^J \DD_r S^{K + 1_i} dL \\
& = (m-1) \sum_p \lambda_p \sum_{|K|=0}^p \sum_{|J|=p-|K|} \sum_r 
\frac{(-1)^r}{J! \, K!} F_{|K|,r} 
d_{J+K} \sum_{q=0}^{\min\{|J|,r\}} G_{|J|,q} \DD_{r-q} S^{J+K + 1_i} dL \\
& = (m-1) \sum_p \lambda_p \sum_{|K|=0}^p \sum_{|J|=p-|K|} \sum_r \sum_{q=0}^{\min\{|J|,r\}} 
\sum_{|N| = r-q} \\
& \hspace{15em} \frac{(-1)^r}{J! \, K! \, N!} F_{|K|,r} G_{|J|,q} 
d_{J+K+N} S^{J+K+N + 1_i} dL \, .
\end{align*}

We now have a complicated coefficient times a sum $d_{J+K+N} S^{J+K+N} S^i dL$, and the task is to replace the combined multi-index $J+K+N$ with a single multi-index, as before. We start by combining $J$ and $K$, setting $s=|K|$, so that we have
\begin{align*}
\lefteqn{\hspace{-2em}(m-1) \sum_p \lambda_p \sum_{s=0}^p \sum_{|K|=s} \sum_{|J|=p-s} \sum_r \sum_{q=0}^{\min\{p-s,r\}} 
\sum_{|N| = r-q} \frac{(-1)^r}{J! \, K! \, N!} \times} \\
& \hspace{15em} F_{s,r} G_{p-s,q} 
d_{J+K+N} S^{J+K+N + 1_i} dL \\
& = (m-1) \sum_p \lambda_p \sum_{s=0}^p  \sum_r \sum_{q=0}^{\min\{p-s,r\}} 
\sum_{|N| = r-q} \frac{(-1)^r}{s! \, (p-s)! \, N!} \times \\
& \hspace{8em} F_{s,r} G_{p-s,q} 
d_N d_{k_1 \cdots k_s k_{s+1} \cdots k_p} S^{k_1 \cdots k_s k_{s+1} \cdots k_p} S^N S^i dL \, .
\end{align*}
Replacing $q$ by $r-q$ gives
\begin{align*}
& (m-1) \sum_p \lambda_p \sum_{s=0}^p  \sum_r \sum_{q=\max\{r-(p-s),0\}}^r \\
& \hspace{8em} \frac{(-1)^r}{s! \, (p-s)! \, q!} F_{s,r} G_{p-s,r-q} 
d_{k_1 \cdots k_p k_{p+1} \cdots k_{p+q}} S^{k_1 \cdots k_p k_{p+1} \cdots k_{p+q}} S^i dL
\end{align*}
and replacing $q$ by $q-p$ gives
\[
(m-1) \sum_p \lambda_p \sum_{s=0}^p \sum_r \sum_{q=\max\{r+s,p\}}^{p+r} 
\frac{(-1)^r}{s! \, (p-s)! \, (q-p)!} F_{s,r} G_{p-s,r-q+p} 
d_{k_1 \cdots k_q} S^{k_1 \cdots k_q} S^i dL \, .
\]

Our task is now to pull the sum over $q$ to the front. The restrictions on $q$ will mean that we have $q \geq r+s$, $q \geq p$ and $q \leq p+r$; we implement these as follows.
\begin{itemize}
\item $q \geq r+s$:\ we already have $s \leq p \leq q$, so this does not constrain $s$. It does, however, constrain $r$, so we set the upper limit of the sum over $r$ to be $q-s$.
\item $q \geq p$:\ we set the upper limit of the sum over $p$ to be $q$.
\item $q \leq p+r$:\ we set the lower limit of the sum over $r$ to be $q-p$.
\end{itemize}
We therefore have
\[
(m-1) \sum_q \sum_{p=0}^q \lambda_p \sum_{s=0}^p \sum_{r=q-p}^{q-s} 
\frac{(-1)^r}{s! \, (p-s)! \, (q-p)!} F_{s,r} G_{p-s,r-q+p} 
d_{k_1 \cdots k_q} S^{k_1 \cdots k_q} S^i dL
\]
and then, replacing $r$ by $r+q-p$, we get
\[
(m-1) \sum_q \sum_{p=0}^q \lambda_p \sum_{s=0}^p \sum_{r=0}^{p-s} 
\frac{(-1)^{r+q-p}}{s! \, (p-s)! \, (q-p)!} F_{s,r+q-p} G_{p-s,r} 
d_{k_1 \cdots k_q} S^{k_1 \cdots k_q} S^i dL
\]
which we shall write as
\[
(m-1) \sum_q H_q \DD_q S^i dL
\]
where
\begin{align*}
H_q & = \sum_{p=0}^q \lambda_p \sum_{s=0}^p \sum_{r=0}^{p-s} 
\frac{(-1)^{r+q-p} \, q!}{s! \, (p-s)! \, (q-p)!} F_{s,r+q-p} G_{p-s,r} \\
& = \sum_{p=0}^q \sum_{s=0}^p \sum_{r=0}^{p-s} 
\frac{(-1)^{r+q} p!q!}{2^{p+1} r!s!(p+2)!(q-p)!(p-s-r)!} \\ 
& \qquad \times \left(\frac{s! (r+q-p)! + (-1)^s (s + r+q-p)!}{(s + r+q-p + 1)!} \right) \\
& = A + B \, .
\end{align*}
We now assert that
\begin{equation}
H_q = \frac{(-1)^q}{2(q+2)} + \frac{(-1)^q (2^{q+1} - 1)  q!}{2^{q+1} (q+2)!} \label{part3}
\end{equation}
and demonstrate this by considering the terms $A$ and $B$ separately.

First, take
\[
A = \sum_{p=0}^q\sum_{s=0}^p\sum_{r=0}^{p-s}
\frac{(-1)^{q+r}p!q!(r+q-p)!}{2^{p+1}(p+2)!r!(q-p)!(p-s-r)!(q+1-(p-s-r))!} \, ;
\]
we shall show that
\[
A = \frac{(-1)^q}{2(q+2)} \, .
\]
Interchange the summations over $r$ and $s$; replace $s$ by $p-r-s$; and reverse the last two summations again:
\[
A=\sum_{p=0}^q\sum_{s=0}^p\frac{(-1)^qp!q!}{2^{p+1}(p+2)!s!(q+1-s)!}
\sum_{r=0}^{p-s}\frac{(-1)^r(r+q-p)!}{r!(q-p)!} \, .
\]
Now consider the sum over $r$:\ it is the coefficient of $x^{q-p}$ in 
\begin{align*}
\sum_{r=0}^{p-s}(-1)^r(1+x)^{r+q-p}
& = (1+x)^{q-p}\sum_{r=0}^{p-s}(-(1+x))^r\\
& = \frac{(1+x)^{q-p}+(-1)^{p-s}(1+x)^{q-s+1}}{2+x} \, .
\end{align*}
The right-hand side is a polynomial of degree $q-s$. Write 
it as $b_0+b_1x+\cdots+b_{q-s}x^{q-s}$:\ then 
\[
(2+x)(b_0+b_1x+\cdots+b_{q-s}x^{q-s})= (1+x)^{q-p}+(-1)^{p-s}(1+x)^{q-s+1} \, .
\]
We want $b_{q-p}$. By comparing coefficients of powers of $x$ we obtain a recurrence relation for the $b$\,s which is easily solved to give
\[
b_{q-p} = \sum_{r=0}^{p-s}(-1)^r2^{p-r-s}\frac{(q+1-s)!}{(q+1-r-s)!r!} \, .
\]
But now we have
\[
A = \sum_{p=0}^q\sum_{s=0}^p\sum_{r=0}^{p-s}
\frac{(-1)^{q+r}p!q!}{2^{r+s+1}(p+2)!r!s!(q+1-r-s)!} \, .
\]
Move the summation over $p$ through to the right:
\[
A=
\sum_{s=0}^q\sum_{r=0}^{q-s}
\frac{(-1)^{q+r}q!}{2^{r+s+1}r!s!(q+1-r-s)!}\sum_{p=r+s}^q\frac{1}{(p+1)(p+2)} \, .
\]
Now
\[
\sum_{p=r+s}^q\frac{1}{(p+1)(p+2)}=
\sum_{p=r+s}^q\left(\frac{1}{(p+1)}-\frac{1}{(p+2)}\right)
=\frac{1}{(r+s+1)}-\frac{1}{(q+2)} \, ,
\]
and thus
\[
A=
\sum_{s=0}^q\sum_{r=0}^{q-s}
\frac{(-1)^{q+r}q!}{2^{r+s+1}r!s!(q+1-r-s)!}
\left(\frac{1}{(r+s+1)}-\frac{1}{(q+2)}\right) \, .
\]
Now set $r+s=t$ in the right-hand summation, and reverse the order of summation to obtain
\[
A = \sum_{t=0}^q\sum_{s=0}^{t}
\frac{(-1)^{q+s+t}q!}{2^{t+1}(t-s)!s!(q+1-t)!}
\left(\frac{1}{(t+1)}-\frac{1}{(q+2)}\right) \, .
\]
But
\[
\sum_{s=0}^{t}\frac{(-1)^s}{(t-s)!s!}=0\quad\mbox{for $t>0$} \, ,
\]
so only the term with $t=0=s$ remains, and therefore
\[
A=\frac{(-1)^{q}q!}{2(q+1)!}\left(1-\frac{1}{(q+2)}\right)
=\frac{(-1)^{q}}{2(q+1)}\left(\frac{q+1}{q+2}\right)
=\frac{(-1)^{q}}{2(q+2)}
\]
as required.

Now consider
\[
B = \sum_{p=0}^q \sum_{s=0}^p \sum_{r=0}^{p-s} 
\frac{(-1)^{r+q+s} p!q!}{2^{p+1} r!s!(p+2)!(q-p)!(p-s-r)!(s + r+q-p + 1)} \, .
\]
Replacing $r$ by $p-r-s$ and interchanging the sums over $s$ and $r$ gives
\[
B = \sum_{p=0}^q \sum_{r=0}^p 
\frac{(-1)^{p+q-r} p!q!}{2^{p+1} (p+2)!(q-p)!r!(q-r+1)}
\sum_{s=0}^{p-r} \frac{1}{(p-r-s)!s!} \, ;
\]
we may then use
\[
\sum_{s=0}^{p-r} \frac{1}{(p-r-s)!s!} = \frac{2^{p-r}}{(p-r)!}
\]
to give
\[
B = \sum_{p=0}^q \sum_{r=0}^p 
\frac{(-1)^{p+q-r} p!q!}{2^{r+1} (p+2)!(q-p)!r!(q-r+1)(p-r)!} \, .
\]
Next, interchange the sums over $p$ and $r$ and then replace $p$ by $p+r$ to give
\[
B = \sum_{r=0}^q 
\frac{(-1)^q q!}{2^{r+1} r!(q-r+1)}
\sum_{p=0}^{q-r} \frac{(-1)^p (p+r)!}{(p+r+2)!p!(q-r-p)!} \, .
\]
The sum over $p$ may be obtained by evaluating the double integral
\[
\int_0^1 \left( \int_0^y x^r (1-x)^{q-r} dx \right) dy
\]
in two different ways. Expanding the inner bracket gives
\begin{align*}
\lefteqn{\hspace{-2em}\sum_{p=0}^{q-r} \frac{(-1)^p (q-r)!}{p! (q-r-p)!} \int_0^1 
\left( \int_0^y x^{p+r} dx \right) dy} \\
& = \sum_{p=0}^{q-r} \frac{(-1)^p (q-r)!}{(p+r+1)(p+r+2)p! (q-r-p)!}
\end{align*}
which is $(q-r)!$ times the sum we want; on the other hand, repeated integration by parts with respect to $x$ gives
\[
\sum_{p=0}^{q-r} \frac{(q-r)!r!}{(q+1-p)!p!} \int_0^1 y^{q+1-p} (1-y)^p dy \, .
\]
But
\[
\int_0^1 y^{q+1-p} (1-y)^p dy = \frac{(q+1-p)!p!}{(q+2)!} \, ,
\]
again using repeated integration by parts, so we have
\[
\sum_{p=0}^{q-r} \frac{(q-r)!r!}{(q+2)!} = \frac{(q-r+1)!r!}{(q+2)!}
\]
for the double integral, and therefore
\[
\sum_{p=0}^{q-r} \frac{(-1)^p (p+r)!}{(p+r+2)!p!(q-r-p)!}
= \frac{(q-r+1)r!}{(q+2)!} \, .
\]
Thus, returning to our original calculation, we have
\[
B = \sum_{r=0}^q \frac{(-1)^q q!}{2^{r+1}(q+2)!} 
= \frac{(-1)^q (2^{q+1} - 1)q!}{2^{q+1}(q+2)!}
\]
as asserted.

We may finally combine equations~(\ref{part1}), (\ref{part2}) and~(\ref{part3}) to give
\begin{align*}
\lefteqn{\hspace{-2em}i_j \left( \sum_{p=0}^\infty \lambda_p \DD_p S^j d\vartheta^i 
- \sum_{p=0}^\infty \lambda_p \DD_p S^i \vartheta^j \right)} \\
& = (m-1) \sum_{p=0}^\infty \frac{(-1)^p \, p!}{2^{p+1} (p+2)!} \DD_p S^i dL + (m-1) \sum_{p=0}^\infty \frac{(-1)^p}{2 \, (p+2)} \DD_p S^i dL \\
& \qquad + \, (m-1) \sum_{p=0}^\infty \left(\frac{(-1)^p}{2(p+2)} + \frac{(-1)^p (2^{p+1} - 1)  p!}{2^{p+1} (p+2)!}\right) \DD_p S^i dL \\
& = (m-1) \sum_{p=0}^\infty \frac{(-1)^p (p! + 2^{p+1} (p+1)! + (2^{p+1}-1)p!)}{2^{p+1} (p+2)!} \DD_p S^i dL \\
& = (m-1) \sum_{p=0}^\infty \frac{(-1)^p p! (p+2)}{(p+2)!} \DD_p S^i dL \\
& = (m-1) \sum_{p=0}^\infty \frac{(-1)^p}{p+1} \DD_p S^i dL \\
& = (m-1) \vartheta^i
\end{align*}
so that $\iT \Theta_2 = (m-1) \Theta_1$ as required. \qed

\textbf{Corollary 5.4}\\
\textit{The fundamental form $\Theta_1$ of a double-integral homogeneous Lagrangian $L$ satisfies the closure condition.} \qed

\section{Discussion}

There are two significant issues raised by the proof of Theorem~5.3 given above, and by the resulting Corollary. The first concerns the restriction, in the present homogeneous formulation, to problems in two independent variables; the second concerns the relationship of these results to variational problems on jet bundles. We consider the second issue first.

It is known that, on a jet bundle, it is always possible to find a Lepage equivalent for any given Lagrangian. There are, however, questions about whether such a Lepage equivalent is unique, and whether --- if not unique --- it nevertheless arises from a natural construction and so transforms appropriately under a change of coordinates. The first of these questions is straightforward, because any two Lepage equivalents must differ by a form which is at least 2-contact, and so the result is unique if, and only if, there is just a single independent variable. The second was answered in~\cite{HK}:\ in the case of several independent variables, there is a natural construction of a Lepage equivalent only in the case of first- and second-order Lagrangians. A similar result must of course hold for Lagrangians defined on bundles of contact elements. But the results for homogeneous Lagrangians are rather different:\ in that context it is always possible to construct an $m$-form which plays the same r\^{o}le as a Lepage equivalent, regardless of the order of the Lagrangian~\cite{CS1}. This $m$-form is, however, projectable to the bundle of contact elements only when $m=1$ or when the Lagrangian is first- or second-order; the result in the case of a first-order Lagrangian is the Carath\'{e}odory form, a natural Lepage equivalent which has been known for many years~\cite{Cth}. 

We are therefore prompted to consider the projectability of the fundamental form of a homogeneous Lagrangian, in the cases where we have demonstrated its existence; of course we need consider only multiple-integral problems. We already know that the fundamental form is projectable for first-order Lagrangians~\cite{CS2}, and we cannot expect such a form to be projectable where the Lagrangian is third-order or higher. There remains, however, the question of second-order Lagrangians. We have demonstrated the existence of a fundamental form for a second-order Lagrangian in two independent variables; if this were projectable then it would give a partial answer to an open question about the existence of fundamental Lepage equivalents for second-order Lagrangians on jet bundles.

In general, however, such a form is not projectable, and a suitable counter-example is given in~\cite{Sau2Order}. The reason for the lack of projectability appears to be that the construction of the fundamental form for two independent variables involves two applications of the homotopy operator $P$, so that the result for a second-order Lagrangian cannot any longer be defined on a third-order bundle; the form constructed then does not survive a reparametrization. This does not, of course, show conclusively that a fundamental Lepage equivalent for a second-order Lagrangian cannot exist; but it is, we believe, evidence that its existence is doubtful. 

We now turn to the issue of extending our results in this paper to the case of homogeneous Lagrangians in more than two independent variables. In principle, the approach would be the same as that adopted above, namely to prove the validity of the remaining stages of the recovery formula. The strategy would also be the same:\ we would express a given vector-valued form $\Theta_r$ in terms of its scalar components $\vartheta^{i_1 \cdots i_r}$ and consider the contractions of these scalar forms with the total derivatives; we would use various commutator formul\ae\ to reduce this to the contractions and Lie derivatives of $\vartheta^{i_1 \cdots i_{r-1}}$ with total derivatives and fundamental vector fields. These in turn would be presumed known from the properties of $\vartheta^{i_1 \cdots i_{r-2}}$, and so on. It is evident that the calculations would be even more complex than those involved in the proof of Theorem~5.3 above. But whereas single-integral variational problems are quite different from multiple-integral ones, we are not aware of any serious difference between problems in two independent variables and problems in three or more, and so we would not expect any inherent obstruction to the construction beyond the complexity just mentioned. Work therefore continues on the project. 
 
\section*{Acknowlegements}

The first author acknowledges the support of grant
no. 201/06/0922 for Global Analysis and its Applications from the 
Czech Science Foundation.

The second author is a Guest Professor at Ghent
University, and a Visiting Senior Research Fellow of King's College,
University of London:\ he is grateful to the
Department of Mathematical Physics and Astronomy at Ghent and the
Department of Mathematics at King's for their hospitality.

\end{document}